\documentclass[12pt]{article}

\usepackage{amsmath,amssymb,amsthm,amscd,a4wide,upref}

\usepackage[enableskew]{youngtab}

\usepackage{hyperref}
\hypersetup{colorlinks,citecolor=blue,filecolor=black,linkcolor=blue,urlcolor=blue}
\usepackage{enumerate}
\usepackage{mathtools}
\usepackage{enumitem}

\usepackage{lmodern}     % set math font to Latin modern math
\usepackage[T1]{fontenc}
\begin{document}

% % % % % % % % % % % % % %
%\newcommand{\tma}{\textcolor{magenta}}
%\newcommand{\tre}{\textcolor{red}}
% % % % % % % % % % % % % %

\newcommand{\ad}{{\rm ad}}
\newcommand{\cri}{{\rm cri}}
\newcommand{\row}{{\rm row}}
\newcommand{\col}{{\rm col}}
\newcommand{\Ann}{{\rm{Ann}\ts}}
\newcommand{\End}{{\rm{End}\ts}}
\newcommand{\Rep}{{\rm{Rep}\ts}}
\newcommand{\Hom}{{\rm{Hom}}}
\newcommand{\Mat}{{\rm{Mat}}}
\newcommand{\ch}{{\rm{ch}\ts}}
\newcommand{\chara}{{\rm{char}\ts}}
\newcommand{\diag}{{\rm diag}}
\newcommand{\st}{{\rm st}}
\newcommand{\non}{\nonumber}
\newcommand{\wt}{\widetilde}
\newcommand{\wh}{\widehat}
\newcommand{\ol}{\overline}
\newcommand{\ot}{\otimes}
\newcommand{\la}{\lambda}
\newcommand{\La}{\Lambda}
\newcommand{\De}{\Delta}
\newcommand{\al}{\alpha}
\newcommand{\be}{\beta}
\newcommand{\ga}{\gamma}
\newcommand{\Ga}{\Gamma}
\newcommand{\ep}{\epsilon}
\newcommand{\ka}{\kappa}
\newcommand{\vk}{\varkappa}
\newcommand{\si}{\sigma}
\newcommand{\vs}{\varsigma}
\newcommand{\vp}{\varphi}
\newcommand{\ta}{\theta}
\newcommand{\de}{\delta}
\newcommand{\ze}{\zeta}
\newcommand{\om}{\omega}
\newcommand{\Om}{\Omega}
\newcommand{\ee}{\epsilon^{}}
\newcommand{\su}{s^{}}
\newcommand{\hra}{\hookrightarrow}
\newcommand{\ve}{\varepsilon}
\newcommand{\pr}{^{\tss\prime}}
\newcommand{\ts}{\,}
\newcommand{\vac}{\mathbf{1}}
\newcommand{\vacu}{|0\rangle}
\newcommand{\di}{\partial}
\newcommand{\qin}{q^{-1}}
\newcommand{\tss}{\hspace{1pt}}
\newcommand{\Sr}{ {\rm S}}
\newcommand{\U}{ {\rm U}}
\newcommand{\BL}{ {\overline L}}
\newcommand{\BE}{ {\overline E}}
\newcommand{\BP}{ {\overline P}}
\newcommand{\AAb}{\mathbb{A}\tss}
\newcommand{\CC}{\mathbb{C}\tss}
\newcommand{\KK}{\mathbb{K}\tss}
\newcommand{\QQ}{\mathbb{Q}\tss}
\newcommand{\SSb}{\mathbb{S}\tss}
\newcommand{\TT}{\mathbb{T}\tss}
\newcommand{\ZZ}{\mathbb{Z}\tss}
\newcommand{\DY}{ {\rm DY}}
\newcommand{\X}{ {\rm X}}
\newcommand{\Y}{ {\rm Y}}
\newcommand{\Z}{{\rm Z}}
\newcommand{\ZX}{{\rm ZX}}
\newcommand{\Ac}{\mathcal{A}}
\newcommand{\Lc}{\mathcal{L}}
\newcommand{\Mc}{\mathcal{M}}
\newcommand{\Pc}{\mathcal{P}}
\newcommand{\Qc}{\mathcal{Q}}
\newcommand{\Rc}{\mathcal{R}}
\newcommand{\Sc}{\mathcal{S}}
\newcommand{\Tc}{\mathcal{T}}
\newcommand{\Bc}{\mathcal{B}}
\newcommand{\Ec}{\mathcal{E}}
\newcommand{\Fc}{\mathcal{F}}
\newcommand{\Gc}{\mathcal{G}}
\newcommand{\Hc}{\mathcal{H}}
\newcommand{\Uc}{\mathcal{U}}
\newcommand{\Vc}{\mathcal{V}}
\newcommand{\Wc}{\mathcal{W}}
\newcommand{\Yc}{\mathcal{Y}}
\newcommand{\Cl}{\mathcal{C}l}
\newcommand{\Ar}{{\rm A}}
\newcommand{\Br}{{\rm B}}
\newcommand{\Ir}{{\rm I}}
\newcommand{\Fr}{{\rm F}}
\newcommand{\Jr}{{\rm J}}
\newcommand{\Or}{{\rm O}}
\newcommand{\GL}{{\rm GL}}
\newcommand{\Spr}{{\rm Sp}}
\newcommand{\Rr}{{\rm R}}
\newcommand{\Zr}{{\rm Z}}
\newcommand{\gl}{\mathfrak{gl}}
\newcommand{\middd}{{\rm mid}}
\newcommand{\ev}{{\rm ev}}
\newcommand{\Pf}{{\rm Pf}}
\newcommand{\Norm}{{\rm Norm\tss}}
\newcommand{\oa}{\mathfrak{o}}
\newcommand{\spa}{\mathfrak{sp}}
\newcommand{\osp}{\mathfrak{osp}}
\newcommand{\f}{\mathfrak{f}}
\newcommand{\g}{\mathfrak{g}}
\newcommand{\h}{\mathfrak h}
\newcommand{\n}{\mathfrak n}
\newcommand{\m}{\mathfrak m}
\newcommand{\z}{\mathfrak{z}}
\newcommand{\Zgot}{\mathfrak{Z}}
\newcommand{\p}{\mathfrak{p}}
\newcommand{\sll}{\mathfrak{sl}}
\newcommand{\agot}{\mathfrak{a}}
\newcommand{\bgot}{\mathfrak{b}}
\newcommand{\qdet}{ {\rm qdet}\ts}
\newcommand{\Ber}{ {\rm Ber}\ts}
\newcommand{\HC}{ {\mathcal HC}}
\newcommand{\cdet}{{\rm cdet}}
\newcommand{\rdet}{{\rm rdet}}
\newcommand{\tr}{ {\rm tr}}
\newcommand{\gr}{ {\rm gr}\ts}
\newcommand{\str}{ {\rm str}}
\newcommand{\loc}{{\rm loc}}
\newcommand{\Gr}{{\rm G}}
\newcommand{\sgn}{ {\rm sgn}\ts}
\newcommand{\sign}{{\rm sgn}}
\newcommand{\ba}{\bar{a}}
\newcommand{\bb}{\bar{b}}
\newcommand{\bi}{\bar{\imath}}
\newcommand{\bj}{\bar{\jmath}}
\newcommand{\bk}{\bar{k}}
\newcommand{\bl}{\bar{l}}
\newcommand{\bp}{\bar{p}}
\newcommand{\hb}{\mathbf{h}}
\newcommand{\Sym}{\mathfrak S}
\newcommand{\fand}{\quad\text{and}\quad}
\newcommand{\Fand}{\qquad\text{and}\qquad}
\newcommand{\For}{\qquad\text{or}\qquad}
\newcommand{\for}{\quad\text{or}\quad}
\newcommand{\grpr}{{\rm gr}^{\tss\prime}\ts}
\newcommand{\degpr}{{\rm deg}^{\tss\prime}\tss}
\newcommand{\bideg}{{\rm bideg}\ts}

\renewcommand{\theequation}{\arabic{section}.\arabic{equation}}

\numberwithin{equation}{section}

\newtheorem{thm}{Theorem}[section]
\newtheorem{lem}[thm]{Lemma}
\newtheorem{prop}[thm]{Proposition}
\newtheorem{cor}[thm]{Corollary}
\newtheorem{conj}[thm]{Conjecture}
\newtheorem*{mthm}{Main Theorem}
\newtheorem*{mthma}{Theorem A}
\newtheorem*{mthmb}{Theorem B}
\newtheorem*{mthmc}{Theorem C}
\newtheorem*{mthmd}{Theorem D}

\theoremstyle{definition}
\newtheorem{defin}[thm]{Definition}

\theoremstyle{remark}
\newtheorem{remark}[thm]{Remark}
\newtheorem{example}[thm]{Example}
\newtheorem{examples}[thm]{Examples}

\newcommand{\bth}{\begin{thm}}
\renewcommand{\eth}{\end{thm}}
\newcommand{\bpr}{\begin{prop}}
\newcommand{\epr}{\end{prop}}
\newcommand{\ble}{\begin{lem}}
\newcommand{\ele}{\end{lem}}
\newcommand{\bco}{\begin{cor}}
\newcommand{\eco}{\end{cor}}
\newcommand{\bde}{\begin{defin}}
\newcommand{\ede}{\end{defin}}
\newcommand{\bex}{\begin{example}}
\newcommand{\eex}{\end{example}}
\newcommand{\bes}{\begin{examples}}
\newcommand{\ees}{\end{examples}}
\newcommand{\bre}{\begin{remark}}
\newcommand{\ere}{\end{remark}}
\newcommand{\bcj}{\begin{conj}}
\newcommand{\ecj}{\end{conj}}

\newcommand{\bal}{\begin{aligned}}
\newcommand{\eal}{\end{aligned}}
\newcommand{\beq}{\begin{equation}}
\newcommand{\eeq}{\end{equation}}
\newcommand{\ben}{\begin{equation*}}
\newcommand{\een}{\end{equation*}}

\newcommand{\bpf}{\begin{proof}}
\newcommand{\epf}{\end{proof}}

\def\beql#1{\begin{equation}\label{#1}}

\newcommand{\Res}{\mathop{\mathrm{Res}}}

\title{\Large\bf Representations of the Yangians associated with\\
Lie superalgebras $\osp(1|2n)$}

\author{A. I. Molev}

\date{} % Start September 2021
\maketitle

%\vspace{4 mm}

\begin{abstract}
We classify the finite-dimensional irreducible representations of
the Yangians associated with the orthosymplectic Lie superalgebras $\osp_{1|2n}$
in terms of the Drinfeld polynomials. The arguments rely on the description
of the representations in the particular case $n=1$ obtained in
our previous work.
\end{abstract}

%\vspace{5 mm}
%%%
%%%{\it Key words:}
%%%

%\newpage

%\tableofcontents
%
%\newpage

\section{Introduction}\label{sec:int}
\setcounter{equation}{0}

The finite-dimensional irreducible representations of
the Yangian $\Y(\g)$ associated with a simple Lie algebra $\g$ were classified by
Drinfeld~\cite{d:nr}. The arguments rely on the work of Tarasov~\cite{t:im}
on the particular case of $\Y(\sll_2)$, where the classification was carried over
in the language of monodromy matrices within the quantum inverse scattering method; see
\cite[Sec.~3.3]{m:yc} for a detailed adapted exposition of these results.
This description of the representations of the Yangian $\Y(\sll_2)$,
along with some other low rank cases,
should also play an essential role
in the classification of the finite-dimensional irreducible representations of
the Yangians associated with simple Lie superalgebras.
One of these cases was considered in our previous work \cite{m:ro}, where
the representations of the Yangian $\Y(\osp_{1|2})$
were described.

These two basic cases turn out to be sufficient to complete the classification in the case of
the Yangians associated with the
orthosymplectic Lie superalgebras $\osp_{1|2n}$. We prove in this paper
that, similar to the classification results of \cite{d:nr},
the finite-dimensional irreducible representations of the Yangian $\Y(\osp_{1|2n})$
are in one-to-one correspondence with the $n$-tuples of monic polynomials
$(P_1(u),\dots,P_n(u))$, and so we call
them the {\em Drinfeld polynomials}.

To describe the results in more detail, recall that
the Yangian
$\Y(\osp_{M|2n})$, as
introduced by Arnaudon {\it et al.\/}~\cite{aacfr:rp}, can be considered
as a quotient
of the extended Yangian $\X(\osp_{M|2n})$ defined via an $RTT$ relation.
A standard argument shows that every finite-dimensional irreducible representation of
$\X(\osp_{M|2n})$ is a highest weight representation. It is isomorphic
to the irreducible quotient $L(\la(u))$ of the Verma module $M(\la(u))$
associated with an $(n+1)$-tuple $\la(u)=(\la_1(u),\dots,\la_{n+1}(u))$ of formal series
$\la_i(u)\in 1+u^{-1}\CC[[u^{-1}]]$. The tuple is called the {\em highest weight}
of the representation.
The key step in the classification is to find
the conditions on the highest weight for
the representation $L(\la(u))$
to be finite-dimensional. The required necessary conditions are derived by
induction from those for the associated
actions of the Yangians $\Y(\gl_2)$ and $\X(\osp_{1|2})$ on the respective
cyclic spans of the highest vector
of $L(\la(u))$. The sufficiency of these conditions is verified by constructing the
{\em fundamental representations} of the Yangian $\X(\osp_{M|2n})$; cf. \cite{amr:rp}, \cite{cp:fr}.
The following is our main result.

\begin{mthm}\label{thm:yclassi}
Every finite-dimensional irreducible representation of the algebra $\X(\osp_{1|2n})$
is isomorphic to $L(\la(u))$ for a certain highest weight $\la(u)$.
The representation $L(\la(u))$
is finite-dimensional if and only if
\beql{ydomire}
\frac{\la_{i+1}(u)}{\la_{i}(u)}=\frac{P_i(u+1)}{P_i(u)},\qquad i=1,\dots,n,
\eeq
for some monic polynomials $P_i(u)$ in $u$. The finite-dimensional irreducible
representations of the Yangian $\Y(\osp_{1|2n})$ are in a one-to-one correspondence with
the $n$-tuples of
monic polynomials $(P_1(u),\dots,P_n(u))$.
\end{mthm}

\section{Definitions and preliminaries}
\label{sec:ns}

For any integer $n\geqslant 1$ introduce the
involution $i\mapsto i\pr=2n-i+2$ on
the set $\{1,2,\dots,2n+1\}$.
Consider the $\ZZ_2$-graded vector space $\CC^{1|2n}$ over $\CC$ with the basis
$e_1,e_2,\dots,e_{2n+1}$, where
the vectors $e_i$ and $e_{i\pr}$ with $i=1,\dots,n$ are odd and
the vector $e_{n+1}$ is even. We set
\ben
\bi=\begin{cases} 1\qquad\text{for}\quad i=1,\dots,n,n',\dots,1',\\
0\qquad\text{for}\quad i=n+1.
\end{cases}
\een
The endomorphism algebra $\End\CC^{1|2n}$ gets a $\ZZ_2$-gradation with
the parity of the matrix unit $e_{ij}$ found by
$\bi+\bj\mod 2$.

We will consider even square matrices with entries in $\ZZ_2$-graded algebras, their
$(i,j)$ entries will have the parity $\bi+\bj\mod 2$.
The algebra of
even matrices over a superalgebra $\Ac$ will be identified with the tensor product algebra
$\End\CC^{1|2n}\ot\Ac$, so that a matrix $A=[a_{ij}]$ is regarded as the element
\ben
A=\sum_{i,j=1}^{2n+1}e_{ij}\ot a_{ij}(-1)^{\bi\tss\bj+\bj}\in \End\CC^{1|2n}\ot\Ac.
\een
We will use the involutive matrix {\em super-transposition} $t$ defined by
$(A^t)_{ij}=A_{j'i'}(-1)^{\bi\bj+\bj}\tss\ta_i\ta_j$,
where we set
\ben
\ta_i=\begin{cases} \phantom{-}1\qquad\text{for}\quad i=1,\dots,n+1,\\
-1\qquad\text{for}\quad i=n+2,\dots,2n+1.
\end{cases}
\een
This super-transposition is associated with the bilinear form on the space $\CC^{1|2n}$
defined by the anti-diagonal matrix $G=[\de_{ij'}\tss\ta_i]$.
We will also regard $t$ as the linear map
\beql{suptra}
t:\End\CC^{1|2n}\to \End\CC^{1|2n}, \qquad
e_{ij}\mapsto e_{j'i'}(-1)^{\bi\bj+\bi}\tss\ta_i\ta_j.
\eeq
In the case of multiple tensor products of the endomorphism algebras,
we will indicate by $t_a$ the map \eqref{suptra}
acting on the $a$-th copy of $\End\CC^{1|2n}$.

A standard basis of the general linear Lie superalgebra $\gl_{1|2n}$ is formed by elements $E_{ij}$
of the parity $\bi+\bj\mod 2$ for $1\leqslant i,j\leqslant 2n+1$ with the commutation relations
\ben
[E_{ij},E_{kl}]
=\de_{kj}\ts E_{i\tss l}-\de_{i\tss l}\ts E_{kj}(-1)^{(\bi+\bj)(\bk+\bl)}.
\een
We will regard the orthosymplectic Lie superalgebra $\osp_{1|2n}$
associated with the bilinear from defined by $G$ as the subalgebra
of $\gl_{1|2n}$ spanned by the elements
\ben
F_{ij}=E_{ij}-E_{j'i'}(-1)^{\bi\tss\bj+\bi}\ts\ta_i\ta_j.
\een
Introduce the permutation operator $P$ by
\ben
P=\sum_{i,j=1}^{2n+1} e_{ij}\ot e_{ji}(-1)^{\bj}\in \End\CC^{1|2n}\ot\End\CC^{1|2n}
\een
and set
\ben
Q=P^{\tss t_1}=P^{\tss t_2}=\sum_{i,j=1}^{2n+1} e_{ij}\ot e_{i'j'}(-1)^{\bi\bj}\ts\ta_i\ta_j
\in \End\CC^{1|2n}\ot\End\CC^{1|2n}.
\een
The $R$-{\em matrix} associated with $\osp_{1|2n}$ is the
rational function in $u$ given by
\ben
R(u)=1-\frac{P}{u}+\frac{Q}{u-\ka},\qquad \ka=-n-1/2.
\een
This is a super-version of the $R$-matrix
originally found in \cite{zz:rf}.
Following \cite{aacfr:rp}, we
define the {\it extended Yangian\/}
$\X(\osp_{1|2n})$
as a $\ZZ_2$-graded algebra with generators
$t_{ij}^{(r)}$ of parity $\bi+\bj\mod 2$, where $1\leqslant i,j\leqslant 2n+1$ and $r=1,2,\dots$,
satisfying certain quadratic relations. In order to write them down,
introduce the formal series
\beql{tiju}
t_{ij}(u)=\de_{ij}+\sum_{r=1}^{\infty}t_{ij}^{(r)}\ts u^{-r}
\in\X(\osp_{1|2n})[[u^{-1}]]
\eeq
and combine them into the matrix $T(u)=[t_{ij}(u)]$ so that
\ben
T(u)=\sum_{i,j=1}^{2n+1} e_{ij}\ot t_{ij}(u)(-1)^{\bi\tss\bj+\bj}
\in \End\CC^{1|2n}\ot \X(\osp_{1|2n})[[u^{-1}]].
\een
Consider the algebra
$\End\CC^{1|2n}\ot\End\CC^{1|2n}\ot \X(\osp_{1|2n})[[u^{-1}]]$
and introduce its elements $T_1(u)$ and $T_2(u)$ by
\ben
T_1(u)=\sum_{i,j=1}^{2n+1} e_{ij}\ot 1\ot t_{ij}(u)(-1)^{\bi\tss\bj+\bj},\qquad
T_2(u)=\sum_{i,j=1}^{2n+1} 1\ot e_{ij}\ot t_{ij}(u)(-1)^{\bi\tss\bj+\bj}.
\een
The defining relations for the algebra $\X(\osp_{1|2n})$ take
the form of the $RTT$-{\em relation}
\beql{RTT}
R(u-v)\ts T_1(u)\ts T_2(v)=T_2(v)\ts T_1(u)\ts R(u-v).
\eeq
As shown in \cite{aacfr:rp}, the product $T(u)\ts T^{\tss t}(u-\ka)$ is a scalar matrix with
\beql{ttra}
T(u-\ka)\ts T^{\tss t}(u)=c(u)1,
\eeq
where $c(u)$ is a series in $u^{-1}$. All its coefficients belong to
the center $\ZX(\osp_{1|2n})$ of $\X(\osp_{1|2n})$ and generate the center.

The {\em Yangian} $\Y(\osp_{1|2n})$
is defined as the subalgebra of
$\X(\osp_{1|2n})$ which
consists of the elements stable under
the automorphisms
\beql{muf}
t_{ij}(u)\mapsto f(u)\ts t_{ij}(u)
\eeq
for all series
$f(u)\in 1+u^{-1}\CC[[u^{-1}]]$.
We have the tensor product decomposition
\beql{tensordecom}
\X(\osp_{1|2n})=\ZX(\osp_{1|2n})\ot \Y(\osp_{1|2n}).
\eeq
The Yangian $\Y(\osp_{1|2n})$ can be equivalently defined as the quotient
of $\X(\osp_{1|2n})$
by the relation
\ben
T(u-\ka)\ts T^{\tss t}(u)=1.
\een

We will also use a more explicit form of the defining relations \eqref{RTT}
written in terms of the series \eqref{tiju} as follows:
\begin{align}
[\tss t_{ij}(u),t_{kl}(v)]&=\frac{1}{u-v}
\big(t_{kj}(u)\ts t_{il}(v)-t_{kj}(v)\ts t_{il}(u)\big)
(-1)^{\bi\tss\bj+\bi\tss\bk+\bj\tss\bk}
\non\\
{}&-\frac{1}{u-v-\ka}
\Big(\de_{k i\pr}\sum_{p=1}^{2n+1}\ts t_{pj}(u)\ts t_{p'l}(v)
(-1)^{\bi+\bi\tss\bj+\bj\tss\bp}\ts\ta_i\ta_p
\label{defrel}\\
&\qquad\qquad\qquad
{}-\de_{l j\pr}\sum_{p=1}^{2n+1}\ts t_{k\tss p'}(v)\ts t_{ip}(u)
(-1)^{\bj+\bp+\bi\tss\bk+\bj\tss\bk+\bi\tss\bp}\ts\ta_j\ta_p\Big).
\non
\end{align}
For any $a\in \CC$ the mapping
\beql{shift}
t_{ij}(u)\mapsto t_{ij}(u+a)
\eeq
defines an automorphism of the algebra $\X(\osp_{1|2n})$.

The universal enveloping algebra $\U(\osp_{1|2n})$ can be regarded as a subalgebra of
$\X(\osp_{1|2n})$ via the embedding
\beql{emb}
F_{ij}\mapsto \frac12\big(t_{ij}^{(1)}-t_{j'i'}^{(1)}(-1)^{\bj+\bi\bj}\ts\ta_i\ta_j\big)(-1)^{\bi}.
\eeq
This fact relies on the Poincar\'e--Birkhoff--Witt theorem for the orthosymplectic Yangian
which was pointed out in \cite{aacfr:rp} and \cite{aacfr:sy}.
It states that the associated graded algebra
for $\Y(\osp_{1|2n})$ is isomorphic to $\U(\osp_{1|2n}[u])$.
A detailed proof of the theorem can be given by
extending the arguments of \cite[Sec.~3]{amr:rp} to the super case
with the use of the vector representation recalled below in \eqref{vectre}.

The extended Yangian $\X(\osp_{1|2n})$ is a Hopf algebra with the coproduct
defined by
\beql{Delta}
\De: t_{ij}(u)\mapsto \sum_{k=1}^{2n+1} t_{ik}(u)\ot t_{kj}(u).
\eeq
For the image of the series $c(u)$ we have $\De:c(u)\mapsto c(u)\ot c(u)$ and so the Yangian
$\Y(\osp_{1|2n})$ inherits the Hopf algebra structure from $\X(\osp_{1|2n})$.

\section{Highest weight representations}
\label{sec:hw}

We will start by deriving a general reduction property for representations of the
extended Yangians $\X(\osp_{1|2n})$ analogous to \cite[Lemma~5.13]{amr:rp}.
For an $\X(\osp_{1|2n})$-module $V$ set
\beql{vplus}
V^+=\{\eta\in V\ |\ t_{1j}(u)\ts\eta=0\quad\text{for}\quad j>1\fand
t_{i\tss 1'}(u)\ts\eta=0\quad\text{for}\quad i<1'\}.
\eeq

\bpr\label{prop:vplus}
The subspace $V^+$ is stable under the action of the operators $t_{ij}(u)$
subject to $2\leqslant i,j\leqslant 2n$.
Moreover, the assignment $\bar t_{ij}(u)\mapsto t_{i+1,j+1}(u)$
for $1\leqslant i,j\leqslant 2n-1$ defines
a representation of the algebra $\X(\osp_{1|2n-2})$ on $V^+$, where the
$\bar t_{ij}(u)$ denote the respective generating series for $\X(\osp_{1|2n-2})$.
\epr

\bpf
Suppose that $2\leqslant k,l\leqslant 2n$ and
$j>1$. For any $\eta\in V^+$ apply \eqref{defrel} to get
\ben
t_{1j}(u)\tss t_{kl}(u)\tss\eta=\frac{1}{u-v-\ka}
\ts\de_{l j\pr}\tss (-1)^{\bj+\bk+\bj\tss\bk}\tss\ta_j\ts t_{k\tss 1'}(v)\ts t_{11}(u)\tss\eta.
\een
Another application of \eqref{defrel} yields
\ben
t_{k\tss 1'}(v)\ts t_{11}(u)\tss\eta=-[t_{11}(u),t_{k\tss 1'}(v)]\ts\eta=
\frac{1}{u-v-\ka}
\ts t_{k\tss 1'}(v)\ts t_{11}(u)\tss\eta,
\een
implying $t_{1j}(u)\tss t_{kl}(u)\tss\eta=0$. A similar calculation shows that
$t_{i\tss 1'}(u)\tss t_{kl}(u)\tss\eta=0$ for $i<1'$ thus proving the first part of
the proposition.

Now suppose that $2\leqslant i,j,k,l\leqslant 2n$. By
\eqref{defrel} the super-commutator $[\tss t_{ij}(u),t_{kl}(v)]$
of the operators in $V^+$ equals
\begin{align}
{}&\frac{1}{u-v}
\big(t_{kj}(u)\ts t_{il}(v)-t_{kj}(v)\ts t_{il}(u)\big)
(-1)^{\bi\tss\bj+\bi\tss\bk+\bj\tss\bk}
\non\\
{}&-\frac{1}{u-v-\ka}
\Big(\de_{k i\pr}\sum_{p=2}^{2n}\ts t_{pj}(u)\ts t_{p'l}(v)
(-1)^{\bi+\bi\tss\bj+\bj\tss\bp}\ts\ta_i\ta_p
\non\\
&\qquad\qquad\qquad\qquad\qquad
{}-\de_{l j\pr}\sum_{p=2}^{2n}\ts t_{k\tss p'}(v)\ts t_{ip}(u)
(-1)^{\bj+\bp+\bi\tss\bk+\bj\tss\bk+\bi\tss\bp}\ts\ta_j\ta_p\Big)
\non
\end{align}
plus the additional terms
\ben
-\frac{1}{u-v-\ka}
\Big(\de_{k i\pr}\ts t_{1j}(u)\ts t_{1'l}(v)
(-1)^{\bi+\bi\tss\bj+\bj}\ts\ta_i
+\de_{l j\pr}\ts t_{k\tss 1'}(v)\ts t_{i1}(u)
(-1)^{\bj+\bi\tss\bk+\bj\tss\bk+\bi}\ts\ta_j\Big).
\een
To transform these terms, use \eqref{defrel} again to get the relations
\ben
\bal
t_{1j}(u)\ts t_{1'l}(v)&=\frac{1}{u-v-\ka-1}\ts
\sum_{p=2}^{2n}\ts t_{pj}(u)\ts t_{p'l}(v)
(-1)^{\bj+\bj\tss\bp}\ts\ta_p\\
&-\frac{1}{u-v-\ka-1}\ts\de_{l j\pr}
\ts t_{1'1'}(v)t_{11}(u)\ts\ta_j
\eal
\een
and
\ben
\bal
t_{k\tss 1'}(v)\ts t_{i1}(u)=[t_{i1}(u)\ts t_{k\tss 1'}(v)](-1)^{\bi+\bk+\bi\tss\bk}&=
\frac{1}{u-v-\ka-1}\ts \de_{k i\pr}\ts t_{11}(u)\ts t_{1'1'}(v)
(-1)^{\bi}\ts\ta_i\\
{}&-\frac{1}{u-v-\ka-1}\ts
\sum_{p=2}^{2n}\ts t_{k\tss p'}(v)\ts t_{ip}(u)
(-1)^{\bi+\bp+\bi\tss\bp}\ts\ta_p.
\eal
\een
Now combine the expressions together and observe that the actions of
the operators $t_{11}(u)$ and $t_{1'1'}(v)$ in $V^+$ commute.
Taking into account the change of the
value $\ka\mapsto\ka+1$ for the algebra $\X(\osp_{1|2n-2})$,
we find that the formula for the super-commutator $[\tss t_{ij}(u),t_{kl}(v)]$
agrees with the defining relations of $\X(\osp_{1|2n-2})$.
\epf

\bre\label{rem:embth}
The reduction property of Proposition~\ref{prop:vplus} should be related to a super-version
of the embedding theorem for the orthogonal and symplectic Yangians
proven in \cite[Thm~3.1]{jlm:ib}. The arguments of that paper should
apply to the super-case to lead to a Drinfeld-type
presentation of the Yangians $\Y(\osp_{1|2n})$ extending the work \cite{aacfr:sy}.
\qed
\ere

A representation $V$ of the algebra $\X(\osp_{1|2n})$
is called a {\em highest weight representation}
if there exists a nonzero vector
$\xi\in V$ such that $V$ is generated by $\xi$,
\begin{alignat}{2}
t_{ij}(u)\ts\xi&=0 \qquad &&\text{for}
\quad 1\leqslant i<j\leqslant 2n+1, \qquad \text{and}\non\\
t_{ii}(u)\ts\xi&=\la_i(u)\ts\xi \qquad &&\text{for}
\quad i=1,\dots,2n+1,
\label{trianb}
\end{alignat}
for some formal series
\beql{laiu}
\la_i(u)\in 1+u^{-1}\CC[[u^{-1}]].
\eeq
The vector $\xi$ is called the {\em highest vector}
of $V$.

\bpr\label{prop:nontrvm}
The series $\la_i(u)$ associated with a highest weight representation $V$
satisfy
the consistency conditions
\beql{nontrvm}
\la_i(u)\tss \la_{i\pr}(u+n-i+1/2)=\la_{i+1}(u)\tss \la_{(i+1)'}(u+n-i+1/2)
\eeq
for $i=1,\dots,n$.
Moreover, the coefficients of the series $c(u)$ act in the representation
$V$ as the multiplications by scalars
determined by
\ben
c(u)\mapsto \la_1(u)\tss \la_{1'}(u+n+1/2).
\een
\epr

\bpf
To derive the consistency conditions,
we will use the induction on $n$ with the base case $n=1$ already considered in \cite{m:ro}.
Suppose that $n\geqslant 2$ and introduce the subspace $V^+$ by \eqref{vplus}.
The vector $\xi$ belongs to $V^+$, and applying Proposition~\ref{prop:vplus} we find that
the cyclic span
$\X(\osp_{1|2n-2})\ts\xi$ is a highest weight submodule with the highest
weight $(\la_2(u),\dots,\la_{2'}(u))$. By the induction hypothesis,
this implies
conditions \eqref{nontrvm} with $i=2,\dots,n$.
Furthermore, using the defining relations \eqref{defrel}, we get
\ben
t_{12}(u)\ts t_{1'2'}(v)\ts \xi=
\frac{1}{u-v-\kappa}\ts\Big(t_{12}(u)\ts t_{1'2'}(v)-
\la_{1}(u)\ts \la_{1'}(v)+\la_{2}(u)\ts \la_{2'}(v)\Big)\ts \xi
\een
and so
\ben
(u-v-\kappa-1)\ts t_{12}(u)\ts t_{1'2'}(v)\ts \xi
=\big({-}\la_{1}(u)\ts \la_{1'}(v)+\la_{2}(u)\ts \la_{2'}(v)\big)\ts \xi.
\een
Setting $v=u-\kappa-1=u+n-1/2$ we obtain
\eqref{nontrvm} for $i=1$. Finally, the last part of the proposition
is obtained by using the expression for $c(u)$ implied by taking the $(1',1')$ entry in
the matrix relation \eqref{ttra}.
\epf

As Proposition~\ref{prop:nontrvm} shows, the series $\la_i(u)$ in \eqref{trianb}
with $i>n+1$ are uniquely
determined by the first $n+1$ series. The corresponding $(n+1)$-tuple
$\la(u)=(\la_{1}(u),\dots,\la_{n+1}(u))$
is called the {\em highest weight\/} of $V$.

Given an arbitrary $(n+1)$-tuple $\la(u)=(\la_{1}(u),\dots,\la_{n+1}(u))$
of formal series of the form \eqref{laiu}, introduce the series $\la_i(u)$
with $i=n+2,\dots,2n+1$ to satisfy the consistency conditions \eqref{nontrvm}.
Define
the {\em Verma module} $M(\la(u))$ as the quotient of the algebra $\X(\osp_{1|2n})$ by
the left ideal generated by all coefficients of the series $t_{ij}(u)$
with $1\leqslant i<j\leqslant 2n+1$, and $t_{ii}(u)-\la_i(u)$ for
$i=1,\dots,2n+1$. As in \cite[Prop.~5.14]{amr:rp},
the Poincar\'e--Birkhoff--Witt theorem for the algebra $\X(\osp_{1|2n})$
implies that the Verma module $M(\la(u))$
is nonzero, and we denote by $L(\la(u))$ its irreducible quotient.

\bpr\label{prop:fdhw}
Every finite-dimensional irreducible representation of the algebra $\X(\osp_{1|2n})$
is isomorphic to $L(\la(u))$ for a certain highest weight
$\la(u)=(\la_{1}(u),\dots,\la_{n+1}(u))$.
\epr

\bpf
The argument is essentially the same as for the proof of the corresponding counterparts
of the property for the Yangians associated with Lie algebras;
cf. \cite[Thm~5.1]{amr:rp}, \cite[Sec.~3.2]{m:yc}. We online some key steps.

Suppose that $V$ is a finite-dimensional irreducible
representation of the algebra $\X(\osp_{1|2n})$ and
introduce its subspace $V^{\tss0}$ by
\ben
V^{\tss0}=\{\eta\in V\ |\ t_{ij}(u)\ts\eta=0,\qquad 1\leqslant i<j\leqslant 2n+1\}.
\een
First we note that $V^{\tss0}$ is nonzero, which follows by considering the set of weights
of $V$, regarded as an $\osp_{1|2n}$-module defined via the embedding
\eqref{emb}. This set is finite and hence contains a maximal weight
with respect to the standard partial ordering on the set of weights of $V$.
A weight vector with this weight belongs to $V^{\tss0}$.

Furthermore, we show that $V^{\tss0}$ is stable under the action of all
operators $t_{ii}(u)$. This follows by straightforward
calculations similar to those used
in the proof of Proposition~\ref{prop:vplus}, relying on the defining relations
\eqref{defrel}. In a similar way, we verify that all the operators
$t_{ii}(u)$ with $i=1,\dots,2n+1$ form a commuting family of operators on $V^{\tss0}$.
Hence they have a simultaneous eigenvector
$\xi\in V^{\tss 0}$. Since the representation $V$ is irreducible,
the submodule $\X(\osp_{1|2n})\tss\xi$ must coincide with $V$
thus proving that $V$ is a highest weight module.

By considering the $\osp_{1|2n}$-weights of $V$ we can also conclude that
the highest vector $\xi$ of $V$ is determined
uniquely, up to a constant factor.
\epf

Proposition~\ref{prop:fdhw} yields the first part of the
Main Theorem. Our next step is to show that the conditions in the theorem
are necessary for the representation $L(\la(u))$ to be finite-dimensional.
So we now suppose that $\dim L(\la(u))<\infty$ and argue by induction on $n$.
The conditions \eqref{ydomire} in the base case $n=1$ are implied by
the main result of \cite{m:ro}. Suppose further that $n\geqslant 2$.

Recall that the Yangian $\Y(\gl_n)$ for the general linear Lie algebra
$\gl_n$ is defined as a unital associative algebra
with countably many generators $t_{ij}^{(1)\circ},\ t_{ij}^{(2)\circ},\dots$ where
$1\leqslant i,j\leqslant n$,
and the defining relations
\ben
(u-v)\ts [t^{\circ}_{ij}(u),t^{\circ}_{kl}(v)]
=t^{\circ}_{kj}(u)\ts t^{\circ}_{il}(v)-t^{\circ}_{kj}(v)\ts t^{\circ}_{il}(u)
\een
written in terms of the series
\ben
t^{\circ}_{ij}(u) = \delta_{ij} + t^{(1)\circ}_{ij} u^{-1} + t^{(2)\circ}_{ij}u^{-2} +
\dots\in\Y(\gl_n)[[u^{-1}]];
\een
see \cite{m:yc} for a detailed exposition of the algebraic structure
and representation of the Yangians associated with $\gl_n$.
The Yangian $\Y(\gl_n)$ can be regarded as a subalgebra of $\X(\osp_{1|2n})$ via the
embedding
\beql{emdy}
\Y(\gl_n)\hra \X(\osp_{1|2n}),\qquad
t^{\circ}_{ij}(u)\mapsto t_{ij}(-u)\quad\text{for}\quad 1\leqslant i,j\leqslant n.
\eeq
The cyclic span $\Y(\gl_n)\tss\xi\subset L(\la(u))$ is a highest weight module
over $\Y(\gl_n)$. Its highest weight is the $n$-tuple $(\la_1(-u),\dots,\la_n(-u))$.
Since $\dim L(\la(u))<\infty$, the corresponding conditions for finite-dimensional
highest weight representations of $\Y(\gl_n)$ must be satisfied; see
\cite[Sec.~3.4]{m:yc}. This implies conditions \eqref{ydomire} of the Main Theorem
for $i=1,\dots,n-1$.

Furthermore, by Proposition~\ref{prop:vplus}, the subspace
$L(\la(u))^+$ is a module over the extended Yangian
$\X(\osp_{1|2n-2})$. The vector $\xi$ generates a highest weight
$\X(\osp_{1|2n-2})$-module
with the highest weight $(\la_{2}(u),\dots,\la_{n+1}(u))$.
Since this module is finite-dimensional,
conditions \eqref{ydomire} hold for $i=2,\dots,n$ by the induction hypothesis.
This completes the proof of the necessity of the conditions.

Now suppose that conditions \eqref{ydomire} hold and derive that the corresponding
module
$L(\la(u))$ is finite-dimensional.
The $n$-tuple of Drinfeld polynomials $(P_1(u),\dots,P_n(u))$ determines the highest weight
$\la(u)$ up to a simultaneous multiplication of all components
$\la_i(u)$ by a series $f(u)\in 1+u^{-1}\CC[[u^{-1}]]$.
This operation corresponds to twisting the action of
the algebra $\X(\osp_{1|2n})$ on $L(\la(u))$ by the automorphism \eqref{muf}.
Hence, it suffices to prove that
a particular module $L(\la(u))$ corresponding to a given set
of Drinfeld polynomials is finite-dimensional.

Suppose that $L(\la(u))$ and $L(\mu(u))$ are the irreducible
highest weight modules with the highest weights
\ben
\la(u)=\big(\la_{1}(u),\dots,\la_{n+1}(u)\big)\Fand
\mu(u)=\big(\mu_{1}(u),\dots,\mu_{n+1}(u)\big).
\een
By the coproduct rule
\eqref{Delta},
the cyclic span $\X(\osp_{1|2n})(\xi\ot\xi')$ of the tensor product
of the respective highest vectors of $L(\la(u))$ and $L(\mu(u))$
is a highest weight
module
with the highest weight
\ben
\big(\la_{1}(u)\tss\mu_{1}(u),\dots,\la_{n+1}(u)\tss\mu_{n+1}(u)\big).
\een
This observation implies that
the cyclic span corresponds to
the set of Drinfeld polynomials $(P_1(u)\tss Q_1(u),\dots,P_n(u)\tss Q_n(u))$,
where the $P_i(u)$ and $Q_i(u)$ are the Drinfeld polynomials for
$L(\la(u))$ and $L(\mu(u))$, respectively.
Therefore,
we only need to establish the sufficiency of conditions \eqref{ydomire}
for the {\em fundamental representations} of $\X(\osp_{1|2n})$ associated
with the $n$-tuples of Drinfeld polynomials
such that
$P_j(u)=1$ for all $j\ne i$ and $P_i(u)=u+b$ for a certain $i\in\{1,\dots,n\}$
and $b\in\CC$; cf.~\cite{cp:fr}.
Moreover, it is sufficient to take one
particular value of $b\in\CC$; the general case will then follow
by twisting the action of the algebra $\X(\osp_{1|2n})$
in such representations
by automorphisms of the form \eqref{shift}.

Consider the vector representation of $\X(\osp_{1|2n})$ on $\CC^{1|2n}$ defined by
\beql{vectre}
t_{ij}(u)\mapsto \de_{ij}+u^{-1}\tss
e_{ij}(-1)^{\bi}-(u+\ka)^{-1}\tss e_{j'i'}(-1)^{\bi\bj}\ts\ta_i\ta_j.
\eeq
The homomorphism property
follows from \eqref{RTT} by applying the standard transposition to one copy of
$\End\CC^{1|2n}$ in the Yang--Baxter equation satisfied by $R(u)$.
Now use the coproduct \eqref{Delta} and suitable automorphisms \eqref{shift}
to equip
the tensor product space $(\CC^{1|2n})^{\ot k}$
with the action of $\X(\osp_{1|2n})$ by setting
\beql{tijtprei}
t_{ij}(u)\mapsto
\sum_{a_1,\dots,a_{k-1}=1}^{2n+1} t_{ia_1}(u)\ot t_{a_1a_2}(u-1)
\ot\dots\ot t_{a_{k-1}j}(u-k+1),
\eeq
where the generators act in the respective copies of the vector space
$\CC^{1|2n}$ via the rule \eqref{vectre}.
For the values $k=1,\dots,n$ introduce the vectors
\ben
\xi_k=\sum_{\si\in\Sym_k} \sgn\si\cdot
e_{\si(1)}\ot\cdots\ot e_{\si(k)} \in
(\CC^{1|2n})^{\ot k}.
\een
Now verify that each vector $\xi_k$ has the properties
\beql{tijxi}
t_{ij}(u)\ts \xi_k=0\qquad\text{for}\quad 1\leqslant i<j\leqslant n+1
\eeq
and
\beql{tiixi}
t_{ii}(u)\ts \xi_k=\begin{cases}
\dfrac{u-k}{u-k+1}\ts\xi_k \qquad&\text{for}
\quad i=1,\dots, k,\\[0.6em]
\ \xi_k\qquad&\text{for}\quad i=k+1,\dots,n+1.
\end{cases}
\eeq
The expression for the vector $\xi_k$ involves
only tensor products of the basis vectors $e_i$ with $i\leqslant n$. This implies
that for the application of the operators $t_{ij}(u)$ with
$1\leqslant i\leqslant j\leqslant n$
to $\xi_k$ we may restrict
the sum in formula \eqref{tijtprei} to the values $a_p\in\{1,\dots,n\}$.

By using the embedding \eqref{emdy}, we may regard the cyclic span $\Y(\gl_n)\tss\xi_k$
as a $\Y(\gl_n)$-module. Moreover, this module is isomorphic to
$A^{(k)}(\CC^n)^{\ot k}$, where $A^{(k)}$ is the anti-symmetrization operator.
It is well-known that this $\Y(\gl_n)$-module is isomorphic to the evaluation
module $L(1,\dots,1,0,\dots,0)$ (with $k$ ones) twisted by
a shift automorphism $u\mapsto u+k-1$; see e.g. \cite[Sec.~6.5]{m:yc}.
This yields formulas \eqref{tijxi} and \eqref{tiixi}
with $1\leqslant i\leqslant j\leqslant n$. They are easily verified directly
for the remaining generators.

Formulas \eqref{tiixi} show that the corresponding set of Drinfeld polynomials
for the highest weight module $\X(\osp_{1|2n})\tss\xi_k$ has the form
$P_i(u)=1$ for $i\ne k$, while $P_k(u)=u-k$. This completes the proof
of the second part of the Main Theorem concerning conditions \eqref{ydomire}.
The last part is immediate from the decomposition
\eqref{tensordecom}; cf. \cite[Sec.~5.3]{amr:rp}.

%\newpage
\bigskip\bigskip

\small

\noindent
School of Mathematics and Statistics\newline
University of Sydney,
NSW 2006, Australia\newline
alexander.molev@sydney.edu.au


\begin{thebibliography}{99}

\bibitem{aacfr:rp}
{D. Arnaudon, J. Avan, N. Cramp\'e, L. Frappat, E. Ragoucy},
{\it $R$-matrix presentation for super-Yangians $Y({\rm osp}(m\vert 2n))$},
{J. Math. Phys.}  {\bf 44}  (2003), 302--308.

\bibitem{aacfr:sy}
{D. Arnaudon, N. Cramp\'e, L. Frappat, E. Ragoucy},
{\it Super Yangian $\Y(osp(1|2))$ and the universal $R$-matrix
of its quantum double},
Comm. Math. Phys. {\bf 240} (2003), 31--51.

\bibitem{amr:rp}
{D. Arnaudon, A. Molev and E. Ragoucy},
{\it On the $R$-matrix realization of Yangians
and their representations},
Annales Henri Poincar\'e {\bf 7} (2006), 1269--1325.

\bibitem{cp:fr}
{V. Chari and A. Pressley},
{\it Fundamental representations of Yangians and rational $R$-matrices},
{J. Reine Angew. Math.} {\bf 417} (1991), 87--128.

\bibitem{d:nr}
{V. G. Drinfeld},
{\it A new realization of
Yangians and quantized affine algebras}, {Soviet Math. Dokl.}
{\bf 36} (1988), 212--216.

\bibitem{jlm:ib}
{N. Jing, M. Liu and A. Molev},
{\it Isomorphism between the $R$-matrix and Drinfeld
presentations of Yangian in types $B$, $C$ and $D$},
Comm. Math. Phys. {\bf 361} (2018), 827--872.

\bibitem{m:yc}
A. Molev,
{\it Yangians and classical Lie algebras}, Mathematical
Surveys and Monographs, 143. AMS,
Providence, RI, 2007.

\bibitem{m:ro}
A. Molev,
{\it Representations of the orthosymplectic Yangian},
{\tt arXiv:2108.10104}.

\bibitem{t:im}
{V. O. Tarasov},
{\it Irreducible monodromy matrices for the $R$-matrix of the
$XXZ$-model and lattice local quantum Hamiltonians}, {Theor. Math. Phys.}
{\bf 63} (1985),
440--454.

\bibitem{zz:rf}
{A. B. Zamolodchikov and Al. B. Zamolodchikov},
{\it Factorized $S$-matrices in two dimensions as the exact solutions
of certain relativistic quantum field models},
{Ann. Phys.} {\bf 120} (1979), 253--291.

\end{thebibliography}
\end{document}